\documentclass[11pt,twoside,a4paper]{article}
\usepackage{amsfonts}
\usepackage{amssymb}
\usepackage{amsmath}
\usepackage[french, english]{babel}

\begin{document}

\newcommand{\pa}{\partial}
\newcommand{\opa}{\overline\pa}
\newcommand{\ol}{\overline }

\numberwithin{equation}{section}

\newcommand\C{\mathbb{C}}  
\newcommand\R{\mathbb{R}}
\newcommand\Z{\mathbb{Z}}
\newcommand\N{\mathbb{N}}
\newcommand\PP{\mathbb{P}}

{\LARGE \centerline{Non locally trivializable $CR$ line bundles}}
{\LARGE \centerline{over compact Lorentzian $CR$ manifolds}}
\vspace{0.8cm}

\centerline{\textsc {Judith Brinkschulte\footnote{Universit\"at Leipzig, Mathematisches Institut, Augustusplatz 10, D-04109 Leipzig, Germany. 
E-mail: brinkschulte@math.uni-leipzig.de}}
 and C. Denson Hill
\footnote{Department of Mathematics, Stony Brook University, Stony Brook NY 11794, USA. E-mail: dhill@math.stonybrook.edu
}}

\vspace{0.5cm}

\begin{abstract}  We consider compact $CR$ manifolds of arbitrary $CR$ codimension that satisfy certain geometric conditions in terms of their Levi form. Over these compact $CR$ manifolds, we construct a deformation of the trivial $CR$ line bundle over $M$ which is topologically trivial over $M$ but fails to be even locally $CR$ trivializable over any open subset of $M$.
In particular, our results apply to compact Lorentzian $CR$ manifolds of hypersurface type.
\end{abstract}

\begin{otherlanguage}{french}
\begin{abstract}
On consid\`ere une vari\'et\'e $CR$ compacte de codimension $CR$ quelconque qui v\'erifie certaines conditions g\'eom\'etriques en terme de sa forme de Levi. Sur ces vari\'et\'es $CR$ compactes, on construit une d\'eformation  du fibr\'e en droites $CR$ trivial sur $M$ qui est topologiquement trivial sur $M$ mais qui n'admet m\^eme pas de trivialization $CR$ locale sur un ouvert arbitraire de $M$. En particulier, nos r\'esultats s'appliquent au cas de vari\'et\'es $CR$ compactes Lorentziennes  du type hypersurface.
\end{abstract}
\end{otherlanguage}

\section{Introduction}

In geometry, the concept of vector bundles is quite important. 
In both categories: smooth real vector bundles over differentiable
real manifolds, and holomorphic complex vector bundles over complex
manifolds, there is a theorem which says that any such vector bundle
always has a local trivialization. In the first situation, the
transition functions, or matrices, are $\mathcal{C}^\infty$ smooth, and in the
second situation, they are holomorphic. For this reason one usually ignores the intrinsic definition of a vector bundle
and works directly with a local trivialization. It is a somewhat suprising
fact, which we address here, that in the category of smooth CR vector
bundles over a smooth CR manifold, the analogous theorem does not always
hold. Thus one can have a perfectly good intrinsically defined CR vector
bundle for which there might not be the corresponding CR transition
functions, or matrices.\\

Our main results are as follows:\\

\newtheorem{main}{Theorem}[section]
\begin{main}  \label{main} \ \\
Let $M$ be a compact (abstract) Lorentzian $CR$ hypersurface of $CR$ dimension $n$, $n\not= 3$.  Then there exists a family $(L_a)_{a>0}$  of complex $CR$ line bundles $L_a\longrightarrow M$ converging to the trivial $CR$ line bundle $L_0= M\times\mathbb{C}$, as $a$ tends to $0$, such that $L_a$ is
 differentiably trivial over $M$, but $L_a$ is not locally $CR$ trivializable  over any open set $U$ of $M$ for $a > 0$.
\end{main}

An easy example satisfying the hypothesis of the above theorem is Penrose's null twistor  space, which is the 5-dimensional $M$ of type (2,1) in $\mathbb{C}\mathbb{P}^3$ given in homogeneous coordinates by
$$\vert z_0\vert^2 + \vert z_1\vert^2 = \vert z_2\vert^2 + \vert z_3\vert^3.$$

For $CR$ manifolds of $CR$ dimension $n=2$ but arbitrary codimension, we obtain the following result:\\

\newtheorem{higher}[main]{Theorem}
\begin{higher}   \label{higher}   \ \\
Let $M$ be a compact (abstract) $CR$ manifold of type $(2,k)$ which is $1$-pseudoconcave. Then there exists a family $(L_a)_{a>0}$  of complex $CR$ line bundles $L_a\longrightarrow M$ converging to the trivial $CR$ line bundle $L_0= M\times\mathbb{C}$, as $a$ tends to $0$, such that $L_a$ is
 differentiably trivial over $M$, but $L_a$ is not locally $CR$ trivializable  over any open set $U$ of $M$ for $a > 0$.
\end{higher}

An example of a $CR$ manifold of type $(2,2)$ satisfying the hypothesis of the above theorem is the twistor space of the Fubini-Study metric on $\C\PP^2$, which is the 6-dimensional $M\subset\C\PP^2\times\C\PP^2$ defined by the complex equation
$$z_0\ol w_0 + z_1\ol w_1 + z_2\ol w_2 =0,$$
where $(z_0:z_1:z_2)$ are homogeneous coordinates  in the first factor, and
$(w_0:w_1:w_2)$ are homogeneous coordinates in the second factor.\\

More examples of $CR$ submanifolds satisfying the hypothesis of this theorem can be found in section 4 of \cite{HN1}.\\

In both theorems the same results apply in the setting of $CR$ vector bundles of arbitrary rank, because it suffices to deform the vector bundle in only one fiber direction.\\

We would like to point out that results in the opposite direction have been obtained by Webster \cite{W}. Namely he proved that any $CR$ vector bundle over a strictly pseudoconvex hypersurface of $CR$ dimension $n\geq 7$ is locally trivializable.

\section{Definitions}

Throughout this paper an abstract $CR$ manifold of type $(n,k)$ is a triple $(M, HM, J)$, where $M$ is a smooth real manifold of dimension $2n+k$, $HM$ is a subbundle of rank $2n$ of the tangent bundle $TM$, and $J: HM \rightarrow HM$ is a smooth fiber preserving bundle isomorphism with $J^2= -\mathrm{Id}$. We also require that $J$ be formally integrable; i.e. that we have
$$\lbrack T^{0,1}M,T^{0,1}M\rbrack \subset T^{0,1}M\qquad\mathrm{where}$$ 
$$ T^{0,1}M = \lbrace X+ iJX\mid X\in \Gamma(M,HM)\rbrace \subset \Gamma(M,\mathbb{C}TM),$$
with $\Gamma$ denoting smooth sections. The $CR$ dimension of $M$ is $n\geq 1$ and the $CR$ codimension is $k\geq 1$.\\

As usual, we denote by $\opa_M$ the tangential Cauchy-Riemann operator on $M$. For the precise definition, we refer the reader to \cite{HN1}.\\

We denote by $H^o M=\lbrace \xi\in T^\ast M\mid < X,\xi>=0, \forall X\in H_{\pi(\xi)}M\rbrace$ the {\it characteristic conormal bundle} of $M$. Here $\pi: T M \longrightarrow M$ is the natural projection. To each $\xi\in H^o_p M$, we associate the Levi form at $\xi:$
$$\mathcal{L}_p(\xi, X) = \xi(\lbrack J\tilde X, \tilde X\rbrack )= d\tilde\xi(X,JX) \ \mathrm{for} \ X\in H_p M$$
which is Hermitian for the complex structure of $H_p M$ defined by $J$. Here $\tilde \xi$ is a section of $H^o M$ extending $\xi$ and $\tilde X$ a section of $HM$ extending $X$. \\

A {\it Lorentzian hypersurface}  is a $CR$ manifold $M$ of type $(n,1)$ such that its Levi form has one negative and $n-1$ positive eigenvalues.\\

$M$ is called {\it $1$-pseudoconcave} if for every $p\in M$ and every characteristic conormal direction $\xi$ the Levi form $\mathcal{L}_p(\xi,\cdot)$ has at least one positive and one negative eigenvalue.\\

For the definition of $CR$ vector bundles, we follow \cite{HN2}. Namely,
a {\it complex $CR$ vector bundle of rank $r$} over an abstract $CR$ manifold $M$ of type $(n,k)$ is a smooth complex vector bundle  $E\overset{\pi}\longrightarrow M$ of rank $r$ such that
\begin{enumerate}
\item[(i)] $E$ has the structure of a smooth abstract $CR$ manifold of type $(n+r,k)$,
\item[(ii)] $\pi: E\longrightarrow M$ is a $CR$ submersion,
\item [(iii)] $E\oplus E\ni (\xi_1,\xi_2)\rightarrow \xi_1 +\xi_2\in E$ and $\mathbb{C}\times E \ni (\lambda,\xi)\rightarrow \lambda\xi\in E$ are $CR$ maps.
\end{enumerate}

Let $E\overset{\pi}\longrightarrow M$ be a complex $CR$ vector bundle of rank $r$. A $CR$ trivialization of $E$ over an open set $U$ of $M$ is a $CR$ equivalence between $E_{\mid U}$ and the trivial bundle $U\times\mathbb{C}^r$. We say that $E$ is {\it locally trivializable} over $M$ iff each point of $M$ has an open neighborhood $U$ such that $E_{\mid U}$ is $CR$ equivalent to $U\times\mathbb{C}^r$. \\

When $E$ is locally $CR$ trivializable, one can arrange an open covering $\lbrace U_\alpha\rbrace$ of $M$, and obtain transition functions, or matrices, $\lbrace g_{\alpha\beta}\rbrace$ which are $CR$. Hence in that case the situation is analogous to the case of holomorphic vector bundles over complex manifolds.\\

Note that if $E$ is locally $CR$ embeddable, then $E$ is also locally $CR$ trivializable, by the implicit function theorem. But the local $CR$ embeddability of the base $M$ does not suffice to make $E$ locally $CR$ trivializable, which is the main point of this paper.\\

Another equivalent definition of $CR$ vector bundles, which is used in \cite{W}, can be given in terms of a connection $\opa_E$  on the complex vector bundle $E$; that is, we have a linear differential operator of order one,
$$\opa_E: \mathcal{C}^\infty_{0,q}(M,E)\longrightarrow \mathcal{C}^\infty_{0,q+1}(M,E)$$
satisfying the Leibnitz rule  and $\opa_E \circ\opa_E=0$. $E$ is locally trivializable over $U\subset M$ iff there exists a nonvanishing section $s$ of $E$ over $U$ satisfying $\opa_E s=0$. For more details, we refer the reader to \cite{HN2}.\\

\section{Construction of a nontrivial global cohomology class}

The key point in the construction of the $CR$ line bundles  is the existence of  a global $\opa_M$-closed $(0,1)$ form on $M$ which is not $\opa_M$-exact on any open subset of $M$. First we consider  a single point $p\in M$, and show the existence of a global $\opa_M$-closed form which is not $\opa_M$-exact on any neighborhood of that point. The existence of such a form follows by similar arguments as in \cite{BH} and \cite{BHN}. Therefore, in the proof of the following Theorem, we sketch only the most important ingredients.\\

\newtheorem{onepoint}{Theorem}[section]
\begin{onepoint}   \label{onepoint}   \ \\
Let $M$ be a compact (abstract) $CR$ manifold of type $(n,k)$. Assume that $H^{0,2}(M)$ is Hausdorff and that there exists a point $p$ on $M$ at which there exists a characteristic conormal direction $\xi$ such that $\mathcal{L}_p(\xi,\cdot)$ has one negative and $n-1$ positive eigenvalues. Then there exists a smooth $(0,1)$-form $f$ on $M$ satisfying $\opa_M f=0$ on $M$ such that $f$ is not $\opa_M$-exact on any neighborhood of  $p$ in $M$.
\end{onepoint}

{\it Proof.} By contradiction, we assume that any smooth $(0,1)$-form $f$ on 
$M$ satisfying $\opa_M f=0$ is  $\opa_M$-exact on some neighborhood  $U_f$ of $p$ in $M$.
Using the functional analytic arguments of \cite{AFN} this implies that there exists an open neighborhood $U$ of $p$ in $M$, independent of $f$,  such that any smooth $(0,1)$-form $f$ on 
$M$ satisfying $\opa_M f=0$ on $M$ is  $\opa_M$-exact on $U$. Moreover, by the open mapping theorem for Fr\'echet spaces, we also get an a priori estimate of the following form:
For a fixed compact $K\subset U$ there exists an integer $m_1$ and a constant $C_1 > 0$ such that the solution $u$ to $\opa_M u = f$ on $U$ can be chosen to satisfy
$$\Vert u\Vert_{K,0} \leq C_1 \Vert f\Vert_{m_1}.$$
Here $\Vert\cdot\Vert_m$ denotes the usual $\mathcal{C}^m$ norm on $\mathcal{C}^\infty_{\cdot,\cdot}(M)$, and $\Vert\cdot\Vert_{K,0}$ denotes the usual $\mathcal{C}^0$ norm on $K$.\\

But this implies that we have
\begin{equation} \label{apriori}
\vert\int_K f\wedge g\vert \lesssim \Vert f\Vert_{m_1}\cdot \Vert\opa_M g\Vert_{K,0}
\end{equation}
for all $f\in \mathcal{C}^\infty_{0,1}(M)\cap\mathrm{Ker}\opa_M$ and $g\in\mathcal{C}_{n+k,n-1}^\infty(M)$ with $\mathrm{supp}\ g\subset K$.\\

Using the geometric assumption on $M$ at $p$, namely that  there exists a characteristic conormal direction $\xi$ such that $\mathcal{L}_p(\xi,\cdot)$ has one negative and $n-1$ positive eigenvalues, one can construct forms $f_\lambda\in\mathcal{C}^\infty_{0,1}(M)$, $g_\lambda\in\mathcal{C}^\infty_{n+k,n-1}(M)$ with support in $K$ such that $\opa_M f_\lambda$ and $\opa_M g_\lambda$ are rapidly decreasing with respect to $\lambda$ in the topology of $\mathcal{C}^\infty_{\cdot,\cdot}(M)$ as $\lambda\rightarrow\infty$. On the other hand we have
\begin{equation}  \label{integral}
\vert\int_K f_\lambda\wedge g_\lambda\vert \geq c \lambda^{-n-\frac{k}{2}}
\end{equation}
for some constant $c > 0$. For the details of the construction of these forms, we refer the reader to section 5 of \cite{BHN}.\\

By assumption $H^{0,2}(M)$ is Hausdorff. But this implies that we can solve the equation $\opa_M u_\lambda = \opa_M f_\lambda$ with an estimate
\begin{equation}  \label{estimate}
\Vert u_\lambda\Vert_{m_1} \leq C_2 \Vert\opa_M f_\lambda\Vert_{m_2},
\end{equation}
where $C_2$ is a positive constant and $m_2$ is an integer. Hence $\Vert u_\lambda\Vert_{m_1}$ is rapidly decreasing with respect to $\lambda$. Defining $\tilde f_\lambda = f_\lambda - u_\lambda$, we obtain a global smooth $\opa_M$-closed $(0,1)$-form on $M$.\\

To get a contradiction to our assumption, we use the estimate (\ref{apriori}) with $\tilde{f}_\lambda$ and $g_\lambda$. Namely, (\ref{apriori}) implies that $\int_M \tilde f_\lambda\wedge g_\lambda$ is rapidly decreasing with respect to $\lambda$, whereas (\ref{integral}) implies that 
 $\vert\int_M \tilde f_\lambda\wedge g_\lambda\vert$ is bounded from below by a  polynomial in $1/\lambda$. For the details of these estimates, we again refer the reader to \cite{BHN}. This contradiction proves the Theorem. \hfill$\square$\\

\newtheorem{nontrivial}[onepoint]{Theorem}
\begin{nontrivial}   \label{nontrivial}   \ \\
Let $M$ be a compact  (abstract) $CR$ manifold of type $(n,k)$. Assume that $H^{0,2}(M)$ is Hausdorff and that  for each $p\in M$ there exists a characteristic conormal direction $\xi$ such that $\mathcal{L}_p(\xi,\cdot)$ has one negative and $n-1$ positive eigenvalues. Then there exists a smooth $(0,1)$-form $\omega$ on $M$ satisfying $\opa_M\omega=0$ on $M$ such that $\omega$ is not $\opa_M$-exact on any neighborhood of any point $p\in M$.
\end{nontrivial}

{\it Proof:} We choose a countable dense set of points on $M$, say $S = (p_j)_{j\in\mathbb{N}}$.
Assume by contradiction that for every $\opa_M$-closed form $f$ on $M$ there exists $j\in\mathbb{N}$ such that $f$ is $\opa_M$-exact on a neighborhood $U_j$ of $p_j$. \\
 
For fixed $j$, let $(U_m^j)_{m\in\N}$ denote a fundamental sequence of open neighborhoods of $p_j$ in $M$. To abbreviate notations, we set $\mathcal{Z}(M)= \mathcal{C}^\infty_{0,1}(M)\cap \mathrm{Ker}\opa_M$.\\
Now for each $m,j\in\mathbb{N}$ we define
$$G_m^j = \lbrace (f, u_j)\in \mathcal{Z}(M) \times \mathcal{C}^\infty(U_m^j) \mid f = \opa_M u_j\ \mathrm{on}\ U_m^j\rbrace.$$

This, as a closed subspace of $\mathcal{Z}(M) \times \mathcal{C}^\infty(U_m^j)$, is also a Fr\'echet space. Let $\pi_m^j: G_m^j \longrightarrow
\mathcal{Z}(M)$ be the natural projection. Our assumption implies that 
$$\mathcal{Z}(M) = \bigcup_{m,j\in\mathbb{N}}\pi_m^j (G_m^j).$$

By Baire's category theorem one of the spaces, $\pi_\ell^s(G_\ell^s)$, must be of second category. Then, by the Banach open mapping theorem, the linear continuous map
$$\pi_\ell^s: G_\ell^s\longrightarrow \mathcal{Z}(M)$$
is surjective. But this contradicts Theorem \ref{onepoint}. \hfill$\square$\\

\section{Construction of $CR$ line bundles}

\newtheorem{linebundle}{Theorem}[section]
\begin{linebundle} \ \\   \label{linebundle}
Let $M$ be an (abstract) $CR$ manifold of type $(n,k)$,  and assume that there exists a smooth $(0,1)$-form $\omega$ on $M$ satisfying $\opa_M\omega=0$ on $M$ such that $\omega$ is not $\opa_M$-exact on any neighborhood of any point $p\in M$. Then there exists a family $(L_a)_{a>0}$  of complex $CR$ line bundles $L_a\longrightarrow M$ converging to the trivial $CR$ line bundle $L_0= M\times\mathbb{C}$, as $a$ tends to $0$, such that $L_a$ is
 differentiably trivial over $M$, but $L_a$ is not locally $CR$ trivializable  over any open set $U$ of $M$ for $a > 0$.
\end{linebundle}

{\it Proof.} The arguments in this section follow \cite{H}. \\

 On the differentiably trivial complex  line bundle $L_a = M\times \mathbb{C}_{z}$, we consider the $CR$ structure whose $T^{0,1}L_a$ is defined as follows:
Let $U$ be an open set of $M$ such that $T^{0,1}M$ is spanned over $U$ by $\ol L_1,\ldots,\ol L_n$. We define $T^{0,1}L_a$ to be 
spanned over $U\times \mathbb{C}_{z}$ by the basis
\begin{equation}  \label{definition}
 \left\{ 
\begin{aligned}
 \ol X_0 & = &\frac{\pa}{\pa\ol z} \hspace{4cm}\\
\ol X_j & = & \ol L_j + a\omega(\ol L_j) \frac{\pa}{\pa z},\ j = 1,\ldots, n
\end{aligned} \right.
\end{equation}
This gives well defined $CR$ structure on $L_a$ . 
 The $CR$ line bundle $L_a$ is differentiably trivial over $M$ and converges to the trivial $CR$ line bundle over $M$ as $a$ tends to zero. 
The associated connection of $L_a$ is defined as follows: Since $L_a$ is smoothly trivial, every $s\in\mathcal{C}^\infty_{0,q}(M,L_a)$ is globally defined by a form $\sigma\in\mathcal{C}^\infty_{0,q}(M)$ and a smooth frame $e\in\mathcal{C}^\infty(M,\mathbb{C})$, $s= \sigma\otimes e$. Then $\opa_{L_a}s = (\opa_M\sigma +\omega\wedge\sigma)\otimes e$. Since $\opa_M\omega =0$, this connection satisfies the integrability condition $\opa_{L_a}\circ\opa_{L_a}=0$.\\

However, for $a\not= 0$, a local $CR$ trivialization of $L_a$ near a point of $M$ forces the existence of a local smooth solution $u$ of $\opa_M u = \omega$.

Indeed, a local nonvanishing $\opa_{L_a}$-closed section means  that we have a local nonvanishing smooth function $\sigma$ on some nonempty open subset of $M$ satisfying $\opa_M\sigma + \omega\sigma =0$. After shrinking $U$, we can assume that $u = -\log\sigma$ is well defined on $U$. But $u$ satisfies $\opa_M u = -\frac{\opa_M\sigma}{\sigma}=\omega$, which contradicts the assumption on $\omega$.
Therefore for $a\not=0$, $L_a$ is not locally $CR$ trivializable over any open subset of $M$. \hfill$\square$\\

\section{Proof of the main theorems}

{\it Proof of Theorem \ref{main}.} Since $M$ is a Lorentzian hypersurface, it satisfies the classical condition $Y(q)$ for $q\not=1,n-1$. If $n\not=3$, then $M$ satisfies in particular the condition $Y(2)$. Hence the $\opa_M$-complex is $\frac{1}{2}$-subelliptic in degree $(0,2)$ (see \cite{FK}). It follows that $H^{0,2}(M)$ is finite dimensional, thus Hausdorff. But then the statement of the theorem follows by combining Theorems \ref{nontrivial} and \ref{linebundle}. \hfill$\square$\\

{\it Proof of Theorem \ref{higher}.}  Since $M$ is $1$-pseudoconcave, it follows from the $\varepsilon$-subelliptic estimates proved in \cite{HN1} that the top-degree cohomology group $H^{0,n}(M)$ is finite dimensional. But since $n=2$ we thus have that $H^{0,2}(M)$ is Hausdorff. Therefore we may again conclude by applying Theorems \ref{nontrivial} and \ref{linebundle}. \hfill$\square$\\

{\it Remark.} We conjecture that Theorem \ref{main} also holds for $n=3$. Indeed, if we start with a $CR$ embedded hypersurface, then we can find many $\opa_M$-closed $(0,1)$-forms defined on a neighborhood of a given point which are not $\opa_M$-exact on any neighborhood of that point (see \cite{AFN}). However it seems to be an open problem to show that $H^{0,2}(M)$ is Hausdorff for a Lorentzian hypersurface of $CR$ dimension $3$. Therefore we cannot construct a {\it global} $\opa_M$-closed $(0,1)$ form on $M$ which is not $\opa_M$-exact on $M$.

\vspace{1cm}

{\bf Acknowledements.} The first author was supported by Deutsche Forschungsgemeinschaft (DFG, German Research Foundation, grant BR 3363/2-1).\\

\end{document}